\documentclass[12pt,leqno]{article}

\usepackage[cp850]{inputenc}
\usepackage{amssymb,amsmath,latexsym,amsthm,amsxtra}

\usepackage{graphicx}

\usepackage[english]{babel}

\usepackage{psfrag}

\numberwithin{equation}{section}

\newtheorem{theo}{Theorem}[section]

\newtheorem{prop}[theo]{Proposition}

\newcommand{\ba}{{\bf {a}}}

\newcommand{\bb}{{\bf {b}}}

\newcommand{\bx}{{\bf {x}}}

\newcommand{\zz}{\mathbb {Z}}

\newcommand{\rr}{\mathbb {R}}

\begin{document}

\title{A graph arising in the Geometry of Numbers}

\author{Wolfgang M. Schmidt \and Leonhard Summerer\thanks{second author supported by FWF grant I 3466-N35}
}

\maketitle

\renewcommand\abstractname{R\'{e}sum\'{e}} 

  \begin{abstract}
La g\'{e}ometrie param\'{e}trique des nombres a permis de visualiser les propri\'{e}t\'{e}s d'approximation simultan\'{e}e d'une collection de nombres r\'{e}els \`{a} travers le graphe combin\'{e} des fonctions de certains minimas successifs. Beaucoup d'in\'{e}galit\'{e}s entre les exposants classiques d'approximation simultan\'{e}e peuvent \^{e}tre d\'{e}duits de ces graphes. En particulier, les graphes dits r\'{e}guliers sont parmis les plus importants, notamment pour les cas extr\^{e}mes de certaines de ces in\'{e}galit\'{e}s. Le but de cet article est de d\'{e}finir et de construire la notion de graphes r\'{e}guliers dans le contexte d'approximation pond\'{e}r\'{e}e.
 \end{abstract}

\renewcommand\abstractname{Abstract} 

\begin{abstract}
  The parametric geometry of numbers has allowed to visualize the simultaneous approximation properties of a collection of real numbers through the combined graph of the related successive minima functions. Several inequalities among classical exponents of simultaneous approximation can be guessed by a study of these graphs; in particular the so called regular graph is of major importance as it provides an extremal case for some of these inequalities. The aim of this paper is to define and construct an analogue of the regular graph in the case of weighted simultaneous approximation.\footnote{
    
    
The authors want to thank the referee for the careful study of the manuscript and his useful comments.}
\end{abstract}

\bigskip

\section{Introduction}

We will first explain how certain graphs arise from Diophantine approximation and the Geometry of Numbers. Yet our construction of such graphs, beginning in Section 2, will not require specific knowledge of these topics.

Diophantine approximation deals with simultaneous approximation to linear forms. Given $n=l+m$ with positive $l,m$, Dirichlet's Theorem in the classical case asserts that for linear forms $f_j(\bx)=\xi_{j1}x_1+\ldots+\xi_{jl}x_l$, $(1\leq j\leq m)$ with real coefficients $\xi_{ji}$ and variables $\bx=(x_1,\ldots,x_l)$, there are non-zero points
$$(x_1,\ldots,x_l,y_1,\ldots,y_m)\in\zz^n\leqno(1.1)$$
with
\begin{eqnarray*}
 |x_i|&\leq&e^{qm}\mbox{ }(1\leq i\leq l),\\
  |f_j(\bx)-y_j|&\leq&e^{-lq}\mbox{ }(1\leq j\leq m).
\end{eqnarray*}                       
The points
$$(x_1,\ldots,x_l,f_1(\bx)-y_1,\ldots,f_m(\bx)-y_m)$$
 with $(1.1)$ form a lattice $\Lambda\subseteq\zz^n$ of covolume $1$, and the box $\mathcal{B}(q)$ of points ${\boldsymbol{\eta}}=(\eta_1,\ldots,\eta_n)$ with
\begin{eqnarray*}
 |\eta_i|&\leq&e^{qm}\mbox{ }(1\leq i\leq l),\\
  |\eta_{l+j}|&\leq&e^{-lq}\mbox{ }(1\leq j\leq m)
\end{eqnarray*}  
has volume $2^n$, so that by Minkowski's first theorem on convex bodies, $\mathcal{B}(q)$ contains a non-zero lattice point, i.e. a point of $\Lambda$.

Recently Diophantine approximation with weights has gained increased attention  (see [1],[3],[10]):  in this more general setup we are given non-negative numbers $\alpha_1,\ldots,\alpha_l,\beta_1,\ldots,\beta_m$, not all zero, with
$$\alpha_1+\ldots+\alpha_l=\beta_1+\ldots+\beta_m.\leqno(1.2)$$
The box $\mathcal{B}(q)$, more precisely $\mathcal{B}_{(\alpha,\beta)}(q)$ now, consisting of the points ${\boldsymbol{\eta}}$ having
\begin{eqnarray*}
 |\eta_i|&\leq&e^{\alpha_i q}\mbox{ }(1\leq i\leq l),\\
  |\eta_{l+j}|&\leq&e^{-\beta_j q}\mbox{ }(1\leq j\leq m)
\end{eqnarray*}  
again has volume $2^n$, hence contains a non-zero point of every lattice $\Lambda$ of covolume $1$. In Minkowski's terminology, the first minimum with respect to $\Lambda$ and $\mathcal{B}(q)$ is $\leq 1$.

Let $\Lambda$ be given and denote Minkowski's successive minima with respect to $\Lambda$ and $\mathcal{B}(q)$ by $\lambda_1(q),\ldots,\lambda_n(q)$. The parametric Geometry of Numbers deals with these minima as functions of the parameter $q\geq 0$. However it is easier to work with their logarithms $L_i(q)=\log\lambda_i(q)$ for $1\leq i\leq n$. A system of functions $L_1,\ldots,L_n$ arising from $\Lambda$ and
$${\boldsymbol{\nu}}:=(\alpha_1,\ldots,\alpha_l,-\beta_1,\ldots,-\beta_m),\leqno(1.3)$$
will be called a $(\Lambda,{\boldsymbol{\nu}})$-system and the union of their graphs will be called a $(\Lambda,{\boldsymbol{\nu}})$-graph. The functions of a $(\Lambda,{\boldsymbol{\nu}})$-system already behave well, but not very well. For instance, by Minkowski's second theorem,
$$-\log n!\leq L_1(q)+\ldots+L_n(q)\leq0,$$
but we wished that the sum was identically zero.

We therefore introduce ${\boldsymbol{\nu}}$-systems as $n$-tuples of functions $P_1(q),\ldots,P_n(q)$ defined for $q\geq0$, which are continuous, satisfy $P_1\leq P_2\leq\ldots\leq P_n$ and $P_1(0)=\ldots=P_n(0)=0$. Moreover, each $P_i$ is piecewise linear, with only finitely many linear pieces in any interval with positive end points, and slopes among $\alpha_1,\ldots,\alpha_l,-\beta_1,\ldots,-\beta_m$. Moreover in every interval where each $P_i$ is linear, the slopes of $P_1,P_2,\ldots,P_n$ will be the above numbers in some order. Hence $P_1(q)+\ldots+P_n(q)=0$ by $(1.2)$.

A ${\boldsymbol{\nu}}$-system will be called \emph{proper} if for each $q>0$ and $i$ with $1\leq i\leq n$ having $P_i(q)<P_{i+1}(q)$, the sum of the slopes of $P_1\ldots,P_i$ to the left of $q$ does not exceed the sum of the slopes of $P_1,\ldots,P_i$ to the right of $q$.

The union of the graphs of $P_1,\ldots,P_n$ will be called a ${\boldsymbol{\nu}}$-graph. A ${\boldsymbol{\nu}}$-system and its graph is \emph{regular} if the graph is invariant under some map ${\boldsymbol{\eta}}\mapsto\tau{\boldsymbol{\eta}}$ with $\tau>1$.

D. Roy in important work [5] showed that in the classical case (i.e. when ${\boldsymbol{\nu}}=(m,\underbrace{-1,\ldots,-1}_{m})$), given any $(\Lambda,{\boldsymbol{\nu}})$-system there is a proper ${\boldsymbol{\nu}}$-system having
$$|L_i(q)-P_i(q)|,\mbox{ }(1\leq i\leq n)\leqno(1.4)$$
bounded independently of $q$. Conversely, given a proper ${\boldsymbol{\nu}}$-system, there is a lattice $\Lambda$ such that $(1.4)$ is bounded. This result allows to have optimal transference inequalities between various exponents of Diophantine approximation. Further progress was made by A. Das, L. Fishman, D. Simmons and M. Urba\'{n}ski in [2] by showing that the above relations between $(\Lambda,{\boldsymbol{\nu}})$-systems and proper ${\boldsymbol{\nu}}$-systems hold for the classical case in general. They used this to provide a variational principle that allows to estimate the Hausdorff dimension of many types of sets determined by Diophantine approximation. Finally, Conjecture $2.3$ of [8] says that the above relationship between $(\Lambda,{\boldsymbol{\nu}})$-systems and proper ${\boldsymbol{\nu}}$-systems holds for weighted Diophantine approximation as well. Whenever this holds, many questions of Diophantine approximation, which play in $\rr^n$, can be reduced to questions on graphs in $\rr^2$. An application regarding Diophantine approximation spectra may be found in [6]. 

In studying ${\boldsymbol{\nu}}$-graphs it will be important to know many examples. In [9] a regular graph for the classical case had been presented which provides an example for a system where the optimal bound for the ratio between ordinary and uniform exponents of Diophantine approximation is attained as was established in [4],[7]. Our goal here will be to construct regular ${\boldsymbol{\nu}}$-graphs in the general, i.e. the weighted, case.

\section{ Construction of regular graphs}

Set $k=$ lcm$(l,m)$. Given ${\boldsymbol{\rho}}=(\rho_1,\ldots,\rho_k)$ with each $\rho_i>1$, we will build a regular graph $\mathcal{G}=\mathcal{G}({\boldsymbol{\nu}},{\boldsymbol{\rho)}}$ depending on the parameters ${\boldsymbol{\nu}}$ and ${\boldsymbol{\rho}}$. Our construction will depend on the ordering within $\{\alpha_1,\ldots,\alpha_l\}$ and $\{\beta_1,\ldots,\beta_m\}$.

Set $\sigma_0=1$, $\sigma_r=\rho_1\cdots\rho_r$ for $1\leq r\leq k$ and $\tau=\sigma_k=\rho_1\cdots\rho_k$. For every $t\in\zz$,
$$\tau^t=\tau^t\sigma_0<\tau^t\sigma_1<\ldots<\tau^t\sigma_k=\tau^{t+1}.$$
The points $\ba_t,\bb_t\in\rr^2$ with $t\in\zz$ of [8] will be replaced by $k$-tuples $(\ba_0^t,\ba_1^t,\ldots,\ba_{k-1}^t)$ and $(\bb_0^t,\bb_1^t,\ldots,\bb_{k-1}^t)$, where for $0\leq r<k$
$$\ba_r^t=\tau^t\sigma_r(1,u_r)\mbox{ and }\bb_r^t=\tau^t\sigma_r(1,v_r)\leqno(2.1)$$
for some numbers $u_0,u_1,\ldots,u_{k-1}$ and $v_0,v_1,\ldots,v_{k-1}$. Observe that for given $r$, the points $\ba_r^t$ will lie on a line $\mathcal{L}_r$ of slope $u_r$ emanating from the origin, and the points $\bb_r^t$ on a line $\mathcal{M}_r$ of slope $v_r$.

When $sk\leq r<(s+1)k$ for some $s\in\zz$, so that $r=sk+h$ with $0\leq h<r$, set $\rho_r=\rho_h$, $u_r=u_h$, $v_r=v_h$ and $\sigma_r=\tau^s\sigma_h$. For instance, when $k\leq r<2k$, we have $\sigma_r=\tau\sigma_h=\rho_1\rho_2\cdots\rho_k\rho_1\rho_2\cdots\rho_h$. We again define $\ba_r^t,\bb_r^t$ by $(2.1)$ and note that for $r=sk+h$ we obtain
$$\ba_r^t=\tau^s\ba_h^t,\mbox{ }\bb_r^t=\tau^s\bb_h^t.\leqno(2.2)$$

We write $\alpha_r=\alpha_i$ if $r\equiv i(\!\!\mod l)$ with $1\leq i\leq l$ and $\beta_r=\beta_j$ if $r\equiv j(\!\!\mod m)$ with $1\leq j\leq m$. Notice that $r\equiv r'(\!\!\mod k)$ implies $r\equiv r'(\!\!\mod l)$ and $r\equiv r'$ $(\!\!\mod m)$, so that we obtain the same pair $i,j$ for $r$ and $r'$. Therefore when $k<lm$, we only need $r$ in $0\leq r<k$. Further denote line segments with end points $\ba,\bb$ by $[\ba,\bb]$ and set for $0\leq r<k$
$$\mathcal{A}_r^t=[\ba_r^t,\bb_{r+l}^t],\mbox{ }\mathcal{B}_r^t=[\bb_r^t,\ba_{r+m}^t].$$

\begin{theo}
  Let $\mathcal{G}$ be the union of ${\bf{0}}$ and of all the line segments $\mathcal{A}_r^t$, $\mathcal{B}_r^t$ with $0\leq r<k$ and $t\in\zz$. There are unique $k$-tuples of numbers $(u_0,u_1,\ldots,u_{k-1})$ and $(v_0,v_1,\ldots,v_{k-1})$ such that each $\mathcal{A}_r^t$ has slope $\alpha_{r+1}$, each $\mathcal{B}_r^t$ has slope $-\beta_{r+1}$ and $\mathcal{G}=\mathcal{G}({\boldsymbol{\nu}},{\boldsymbol{\rho}})$ is a regular ${\boldsymbol{\nu}}$-graph.
\end{theo}


Before we proceed with the proof of Theorem 2.1 we indicate the principle of construction of such a graph in the case $(l,m)=(3,2)$ in the interval $[\tau^t,\tau^{t+1}]$ for some $t\in\zz$. 

\vspace{8mm}

\scalebox{0.48}{

\psfrag{T}{\huge Figure 1}

\psfrag{L}{\Large$\!\!\!\mathcal{L}_0$}
\psfrag{M}{\Large$\!\!\!\!\!\!\mathcal{M}_0$} 
\psfrag{A}{\Large $\mathcal{A}_0^t$}
\psfrag{B}{\Large $\mathcal{B}_3^t$} 
\psfrag{a0}{\Large$\bb_0^t$}
\psfrag{a1}{\Large$\bb_1^t$}
\psfrag{a2}{\Large$\bb_2^t$}
\psfrag{a3}{\Large$\bb_3^t$}
\psfrag{a4}{\Large$\bb_4^t$}
\psfrag{a5}{\Large$\bb_5^t$}
\psfrag{a}{\Large$\bb_0^{t+1}=\tau\bb_0^t$}
\psfrag{bo}{\Large$\ba_0^t$}
\psfrag{b1}{\Large$\ba_1^t$}
\psfrag{b2}{\Large$\ba_2^t$}
\psfrag{b3}{\Large$\ba_3^t$}
\psfrag{b4}{\Large$\ba_4^t$}
\psfrag{b5}{\Large$\ba_5^t$}
\psfrag{b}{\Large$\ba_0^{t+1}=\tau\ba_0^t$}
 
\psfrag{t}{\Large$\tau^t$} 
\psfrag{t1}{\Large$\tau^t\sigma_1$} 

\psfrag{t2}{\Large$\tau^t\sigma_2$} 
\psfrag{t3}{\Large$\tau^t\sigma_3$} 
\psfrag{t4}{\Large$\tau^t\sigma_4$} 
\psfrag{t5}{\Large$\tau^t\sigma_5$}
\psfrag{t6}{\Large$\tau^{t+1}$}
\includegraphics{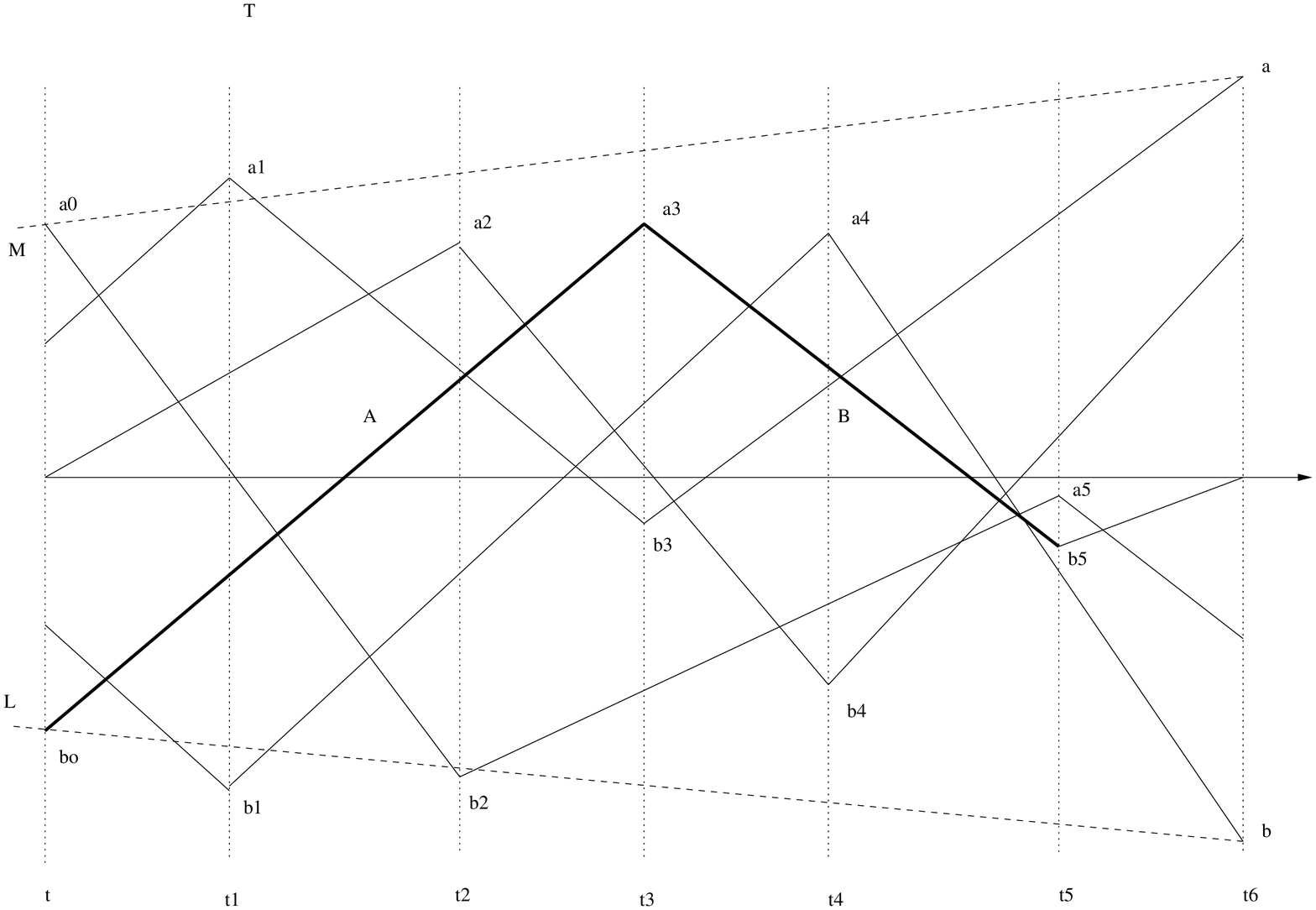}}
\vspace{4mm}

{\bf{Proof}}:
Let us first assume $l,m$ to be relatively prime. Going from $\ba_r^t$ with $0\leq r<k$ to $\bb_{r+l}^t$ via $\mathcal{A}_r^t$, and then via $\mathcal{B}_{r+l}^t$ to $\ba_{r+n}^t$ (look at the bold line segments for $r=0$ in the picture), denoted by
$$\ba_r^t\stackrel{\mathcal{A}_r^t}{\longrightarrow}\bb_{r+l}^t\stackrel{\mathcal{B}_{r+l}^t}{\longrightarrow}\ba_{r+n}^t,$$
(note that $r+l,r+n$ may exceed $k-1$) the ordinate will change by
$$\alpha_{r\!+\!1}(\tau^t\!\sigma_{r+l}\!-\!\tau^t\!\sigma_r)-\beta_{r+l+1}(\tau^t\!\sigma_{r+n}-\tau^t\!\sigma_{r+l})=\tau^t\!\sigma_r(\alpha_{r\!+\!1}(\psi_r^l-1)-\beta_{r+l+1}(\psi_r^n-\psi_r^l)),\leqno(2.3)$$
where $\psi_r^s=\rho_{r+1}\rho_{r+2}\cdots\rho_{r+s}$ (the superscript indicates the number of factors).

On the other hand, the ordinate at $\ba_{r+n}^t-\ba_r^t$ is
$$\tau^t\sigma_{r+n}u_{r+n}-\tau^t\sigma_{r}u_{r}=\tau^t\sigma_r(\psi_r^nu_{r+n}-u_r),$$
so that comparison with $(2.3)$ and division by $\tau^t\sigma_r$ yields
$$u_r-\chi_{r}u_{r+n}=-U_r,\leqno(2.4)$$
with $\chi_r:=\psi_r^n$ and
$$U_r=\alpha_{r+1}(\psi_r^l-1)-\beta_{r+l+1}(\psi_r^n-\psi_r^l).\leqno(2.5)$$
These equations hold for $0\leq r<k$, and yield $k$ linear relations for $u_0,\ldots,u_{k-1}$.

In a similar way, considering the path
$$\bb_r^t\stackrel{\mathcal{B}_r^t}{\longrightarrow}\ba_{r+m}^t\stackrel{\mathcal{A}_{r+m}^t}{\longrightarrow}\bb_{r+n}^t,$$
we obtain
$$v_r-\chi_{r}v_{r+n}=-V_r,\leqno(2.6)$$
with 
$$V_r=-\beta_{r+1}(\psi_r^m-1)+\alpha_{r+m+1}(\psi_r^n-\psi_r^m).\leqno(2.7)$$

As $l,m$ are relatively prime by assumption, so are $n,k$. By $(2.4)$ for $r$, $r+n,r+2n$ we have
\begin{eqnarray*}u_r&=&-U_r+\chi_ru_{r+n}\\
                    &=&-U_r-\chi_rU_{r+n}+\chi_r\chi_{r+n}u_{r+2n}\\
                    &=&-U_r-\chi_rU_{r+n}-\chi_r\chi_{r+n}U_{r+2n}+\chi_r\chi_{r+n}\chi_{r+2n}u_{r+3n}\\
                    &=&\ldots.\end{eqnarray*}
Continuing in this way we obtain after $k$ steps:
$$u_r\!=\!-U_r\!-\!\chi_rU_{r\!+\!n}\!-\!\ldots\!-(\chi_r\chi_{r\!+\!n}\!\cdots\!\chi_{r\!+\!(k-2)n})U_{r\!+\!(k-1)n}\!+\!(\chi_r\chi_{r\!+\!n}\!\cdots\!\chi_{r\!+\!(k-1)n})u_{r\!+\!kn}.\leqno(2.8)$$                 
But $u_{r+kn}=u_r$, and its coefficient in $(2.8)$ is
$$\chi_r\chi_{r+n}\cdots\chi_{r+(k-1)n}=\chi_0\chi_1\cdots\chi_{k-1}=\tau^n$$
since our subscripts are residue classes modulo $k$, and since $n,k$ are coprime, $r,r+n,\ldots,r+(k-1)n$ runs through all these classes, and in view of $\chi_h=\rho_{h+1}\cdots\rho_{h+n}$, each $\rho_i$ occurs in $n$ of the numbers $\chi_0,\ldots,\chi_{k-1}$. Therefore $(2.8)$ is
$$(\tau^n-1)u_r=U_r+\sum_{j=1}^{k-1}(\chi_r\chi_{r+n}\cdots\chi_{r+(j-1)n})U_{r+jn}.\leqno(2.9)$$
The analogous equation for $v_r$, with $U_0,U_1,\ldots,U_{k-1}$ replaced by $V_0,V_1,\ldots,V_{k-1}$ yields
$$(\tau^n-1)v_r=V_r+\sum_{j=1}^{k-1}(\chi_r\chi_{r+n}\cdots\chi_{r+(j-1)n})V_{r+jn}.\leqno(2.9')$$
In combination with $(2.5), (2.7)$ this gives
$$u_r=\frac{1}{\tau^n-1}\sum_{j=0}^{k-1}\left(\alpha_{r+1+jn}(\psi_r^{jn+l}-\psi_r^{jn})-\beta_{r+1+l+jn}(\psi_r^{(j+1)n}-\psi_r^{jn+l})\right),\leqno(2.10)$$
$$v_r=\frac{1}{\tau^n-1}\sum_{j=0}^{k-1}\left(\alpha_{r+1+m+jn}(\psi_r^{(j+1)n}-\psi_r^{jn+m})-\beta_{r+1+jn}(\psi_r^{jn+m}-\psi_r^{jn})\right)\leqno(2.10')$$
and the $u_i,v_i$ are uniquely determined by ${\boldsymbol{\rho}}$ and the ordered set of weights. 

Note that in every interval $[\tau^t_{\sigma_i},\tau^t_{\sigma_{i+1}}]$ the graph 
$\mathcal{G}$ determined by the values of $u_i$ resp. $v_i$, $0\leq i<k$, given by $(2.10)$ resp. $(2.10')$ consists of $n$ line segments having slopes $\alpha_1,\ldots,\alpha_l,-\beta_1,\ldots,-\beta_m$ in some order. Thus $\mathcal{G}$ will be the union of the graphs of $n$ functions $P_1\leq P_2\leq \ldots\leq P_n$, where in each interval $[\tau^t_{\sigma_i},\tau^t_{\sigma_{i+1}}]$ the slopes will be $\alpha_1,\ldots,\alpha_l$ and $-\beta_1,\ldots,-\beta_m$ in some order. Since $\mathcal{A}_r^t$ and $\mathcal{B}_{r+l}^t$ are joined at $b^t_{r+l}$ for $0\leq r<k$, all the functions $P_1,\ldots,P_n$ are continuous. Finally, $\mathcal{G}$ is regular by $(2.2)$, concluding the proof of Theorem 2.1 in the case $(l,m)=1$.\\

In the next section, we will still assume $k,l$ to be relatively prime and complete the proof in Section 4 in the case $(l,m)>1$.

\section{ On proper graphs}

Our regular graph $\mathcal{G}(\boldsymbol{\nu},\boldsymbol{\rho})$ is proper if $v_r\geq u_r$ for $0\leq r<k$, as it is the case for the graph depicted in Figure 1. Moreover we have

\begin{prop}
With
$$\Omega:=\max_{i,j}\{\alpha_i+\beta_j\}\mbox{ ; }\omega:=\min_{i,j}\{\alpha_i+\beta_j\}\mbox{ for }1\leq i\leq l\mbox{ and }1\leq j\leq m.$$
and assuming that
$$\omega>0.\leqno(3.1)$$
we have: $\mathcal{G}(\boldsymbol{\nu},\boldsymbol{\rho})$ is proper provided
$$\frac{\psi_r^n+1}{\psi_r^l+\psi_r^m}\geq\frac{\Omega}{\omega}\mbox{ for any }r\in\{0,\ldots,k-1\}.\leqno(3.2)$$
\end{prop}

Proof: By $(2.5)$, $(2.7)$ $v_r\geq u_r$ for $0\leq r<k$ certainly holds if
$$\alpha_{r+m+1}(\psi_r^n-\psi_r^m)-\alpha_{r+1}(\psi_r^l-1)+\beta_{r+l+1}(\psi_r^n-\psi_r^l)-\beta_{r+1}(\psi_r^m-1)\geq0\leqno(3.3)$$
for each $r$. Now by definition of $\omega$ we have
$$(\alpha_{r+m+1}+\beta_{r+l+1})\psi_r^n+(\alpha_{r+1}+\beta_{r+1})\geq\omega(\psi_r^n+1)$$
and likewise, by definition of $\Omega$:
$$(\alpha_{r+m+1}+\beta_{r+1})\psi_r^m+(\alpha_{r+1}+\beta_{r+l+1})\psi_r^l\leq\Omega(\psi_r^m+\psi_r^l).$$
Thus $(3.3)$ holds for each $r$ provided
$$\omega(\psi_r^n+1)\geq\Omega(\psi_r^l+\psi_r^m)\mbox{ for } r=0,\ldots,k-1.$$
In view of $(3.1)$ this yields exactly $(3.2)$ as claimed.
\vspace{4mm}

Note that from $(3.3)$ it easily follows that for given ${\boldsymbol{\nu}}=({\boldsymbol{\alpha}},{\boldsymbol{\beta}})$, hence given $\Omega/\omega$, the graph $\mathcal{G}({\boldsymbol{\nu}},{\boldsymbol{\rho)}}$ is certainly proper if each $\rho_i$ is sufficiently large. In the case when $m=1$, $\beta_1=l$ we have $\psi_r^n=\tau\rho_{r+1}$, $\psi_r^m=\rho_{r+1}$, $\psi_r^l=\tau$ and $\alpha_{r+m+1}=\alpha_{r+2}$, so that $(3.1)$ becomes
$$\alpha_{r+2}(\tau-1)\rho_{r+1}-\alpha_{r+1}(\tau-1)+l(\tau-1)(\rho_{r+1}-1)\geq0,$$
which is the same as
$$(\alpha_{r+2}+l)\rho_{r+1}-(\alpha_{r+1}+l)\geq0.$$
Replacing $r$ by $r-1$, a simple sufficient condition for being proper is
$$\rho_r\geq\frac{\alpha_r+l}{\alpha_{r+1}+l}\mbox{ for }r=1,\ldots,l.$$

We conclude this section by showing in Figure 2 the example of our regular, proper graph $\mathcal{G}({\boldsymbol{\alpha}},{\boldsymbol{\beta}},{\boldsymbol{\rho}})$ obtained for the parameters 
    \begin{itemize}
    \item $l=3$, $m=2$,\\
    \item $\rho_i=\sqrt[3]{2}$ for $i=1,\ldots,6$, hence $\tau=4$,\\
    \item ${\boldsymbol{\alpha}}=(1/2,1,3/2)$, ${\boldsymbol{\beta}}=(2,1)$.
    \end{itemize}
We restrict the picture to the part with $q$ in the interval $[\tau^0,\tau^1]=[1,4]$. It is plain to see that this graph is proper as any $\bb_r^0=\sigma_r(1,v_r)$ lies above the corresponding $\ba_r^0=\sigma_r(1,u_r)$, so that $v_r>u_r$.
\vspace{8mm}

\scalebox{0.55}{

\psfrag{T}{\huge Figure 2}

\psfrag{A}{\Large slope $v_0$} 
\psfrag{B}{\Large slope $u_0$} 
\psfrag{b}{\Large\:$4$}
\psfrag{t}{\Large$1$} 
\psfrag{t1}{\Large$\sqrt[3]{2}$} 

\psfrag{t2}{\Large$\sqrt[3]{4}$} 
\psfrag{t3}{\Large$2$} 
\psfrag{t4}{\Large$2\sqrt[3]{2}$} 
\psfrag{t5}{\Large$2\sqrt[3]{4}$}
\psfrag{t6}{}
\psfrag{a0}{\Large$\!\bb_0^0$}
\psfrag{a1}{\Large$\bb_1^0$}
\psfrag{a2}{\Large$\!\!\bb_2^0$}
\psfrag{a3}{\Large$\!\!\bb_3^0$}
\psfrag{a4}{\Large$\!\!\bb_4^0$}
\psfrag{a5}{\Large$\bb_5^0$}
\psfrag{a}{\Large$\bb_0^{1}$}
\psfrag{bo}{\Large$\!\ba_0^0$}
\psfrag{b1}{\Large$\ba_1^0$}
\psfrag{b2}{\Large$\!\!\ba_2^0$}
\psfrag{b3}{\Large$\!\ba_3^0$}
\psfrag{b4}{\Large$\!\ba_4^0$}
\psfrag{b5}{\Large$\!\ba_5^0$}
\psfrag{b0}{\Large$\ba_0^{1}$}
 
\includegraphics{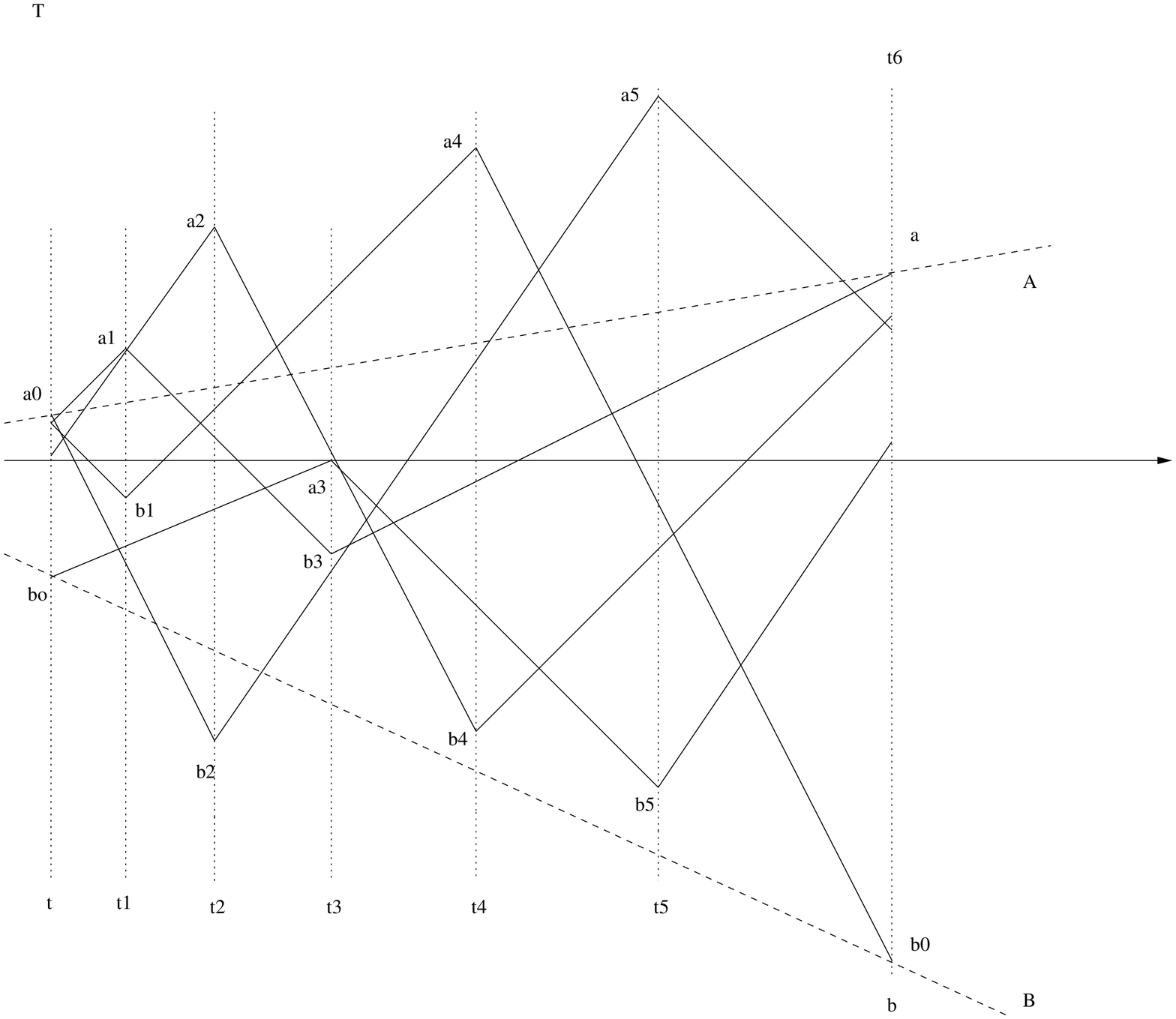}}
\vspace{4mm}

\section{The case when $l,m$ are not coprime}

We now finally assume $k,l$ to have greatest common factor $d>1$, so that $k=lm/d$. We set $n'=n/d$, $k'=k/d=lm/d^2$ and note that $n',k'$ are coprime. Given a residue class $f$ modulo $d$, the line segments
$$\mathcal{A}_{r}^t,\mathcal{B}_r^t\mbox{ with }r\equiv f(\!\!\mod d)\leqno(4.1)$$
have endpoints $\ba_s^t,\bb_w^t$ with $s\equiv w\equiv f (\!\!\mod d)$ and slopes $\alpha_{r+1}^t,-\beta_{r+1}^t$. The union of ${\bf{0}}$ and the line segments $(4.1)$ is a graph $\mathcal{G}^f$, and $\mathcal{G}$ is the union of $\mathcal{G}^1,\ldots,\mathcal{G}^d$. With ${\boldsymbol{\nu}}^f\in\rr^{n'}$ having the components $\alpha_i,-\beta_j$ with $i\equiv j\equiv f(\!\!\mod d)$, the graph $\mathcal{G}^f$ is a ${\boldsymbol{\nu}}^f$-graph, i.e. it is like a ${\boldsymbol{\nu}}$-graph, except that the associated functions $P_1^f,\ldots,P_{n'}^f$ will in every interval where each of them is linear, have the components of ${\boldsymbol{\nu}}^f$ as slopes in some order. Thus when $\gamma^f$ is the sum of these components, we will have $P_1^f(q)+\ldots+P_{n'}^f(q)=\gamma^f q$.

Now let us go back to $(2.4)$. Note that when $h$ runs through the residue classes modulo $k'$, then $f\!+\!dh$ runs through the residue classes modulo $k$ that are congruent to $f (\!\!\mod d)$. Furthermore, $u_h^f:=u_{f+dh}$ will run through the slopes $u_r$ belonging to $\mathcal{G}^f$. With $U_h^f:=U_{f+hd}$, $\chi_h^f:=\chi_{f+hd}$, $(2.4)$ yields
$$u_h^f-\chi_h^fu_{h+n'}^f=-U_h^f.\leqno(4.2)$$

In analogy to $(2.8)$ we have
\begin{eqnarray*}(4.3)u_h^f\!\!&=&\!\!-U_h^f-\chi_h^fU_{h+n'}^f-\chi_h^f\chi_{h+n'}^fU_{h+2n'}^f-\!\ldots\!\\
&& -(\chi_h^f\chi_{h+n'}^f\!\cdots\chi_{h+(k'-2)n'}^f)U_{h+(k'-1)n'}^f+(\chi_h^f\chi_{h+n'}^f\cdots\chi_{h+(k'-1)n'}^f)u_{h+k'n'}^f\end{eqnarray*}
with subscripts modulo $k'$. We have $k'n'\equiv0(\!\!\!\mod k')$ so that $u_{h+k'n'}^f=u_h^f$, with coefficient $\chi_h^f\chi_{h+n'}^f\cdots\chi_{h+(k'-1)n'}^f=\chi_h^f\chi_{h+1}^f\cdots\chi_{h+(k'-1)}^f$, since with $n',k'$ coprime, $0,n',\ldots,(k'-1)n'$ runs through all the residue classes modulo $k'$. Each $\rho_i$ is a factor of $n'$ of the numbers $\chi_0^f,\ldots,\chi_{k'-1}^f$ so that this coefficient is $\tau^{n'}$. Therefore $(4.3)$ gives
$$(\tau^{n'}-1)u_h^f=U_h^f+\sum_{j=1}^{k'-1}(\chi_h^f\chi_{h+n'}^f\cdots\chi_{h+(j-1)n'}^f)U_{h+jn'}^f.\leqno(4.4)$$

An analogous relation holds for $v_h^f=v_{f+dh}$.
\vspace{20mm}

\begin{center}{\Large{\bf References}}\end{center}

\vspace{4mm}

[1] S. Chow, A. Ghosh, L. Guan, A. Marnat, D. Simmons; {\emph{Diophantine transference inequalities: weighted, inhomogeneous, and intermediate exponents}}.
Annali Della Scuola Normale Superiore Di Pisa, arXiv:1808.07184v2.
\vspace{4mm}

[2] T. Das, L. Fishman, D. Simmons, M. Urba\'{n}ski; {\emph{A variational principle in the parametric geometry of numbers}} https://arxiv.org/abs/1901.06602
\vspace{4mm}

[3] O. German; {\emph {Transference theorems for Diophantine approximation with\\ weights}}. Mathematika {\bf 66} No. 2 (2020), p. 325-342
\vspace{4mm}

[4] A. Marnat, N. Moshchevitin; {\emph{An optimal bound for the ratio between ordinary and uniform exponents of Diophantine approximation}} Mathematika {\bf 66} No. 3 (2020), p. 818-854
\vspace{4mm}

[5] D. Roy; {\emph{On Schmidt and Summerer parametric geometry of numbers}}. Ann. of Math.{\bf 182} (2015), p. 739-786
\vspace{4mm}

[6] D. Roy; {\emph{On the topology of Diophantine approximation Spectra}}. Compositio Math. {\bf 153} (2017), p. 1512-1546
\vspace{4mm}

[7] M. Rivard-Cooke; {\emph{Parametric Geometry of Numbers}}. PhD Thesis, University of Ottawa (2019) https://ruor.uottawa.ca/handle/10393/38871
\vspace{4mm}

[8] W.M. Schmidt; {\emph{On Parametric Geometry of Numbers}}. Acta Arithmetica {\bf 195} (2020) P. 383-414
\vspace{4mm}

[9] W. M. Schmidt, L. Summerer; {\emph{Parametric Geometry of Numbers and applications}}. Acta Arithmetica  {\bf 140} No. 1 (2009), p. 67-91 
\vspace{4mm}

[10] L. Summerer; {\emph{A geometric proof of Jarnik's identity in the setting of weighted simultaneous approximation}} https://arxiv.org/abs/1912.04574
\vspace{8mm}

Wolfgang M. Schmidt\\
Department of Mathematics\\
University of Colorado \\
Boulder, CO 80309-0395, USA\\
e-mail: wolfgang.schmidt@colorado.edu 
\vspace{5mm}

Leonhard Summerer\\
Fakult\"at f\"ur Mathematik\\
Universit\"at Wien\\
Oskar-Morgenstern-Platz 1\\
A-1090 Wien, AUSTRIA\\
e-mail: leonhard.summerer@univie.ac.at

\end{document}